\documentclass{article}
\usepackage{ucs}
\usepackage[utf8]{inputenc}

\usepackage{grffile}
\usepackage{ctable}

\usepackage[small]{caption}
\usepackage{graphicx}

\usepackage{amsmath,amsfonts,amssymb}
\usepackage{amsthm}
\usepackage{color}
\usepackage{mathrsfs}
\usepackage{stmaryrd}
\usepackage{fullpage}
\usepackage{ifthen}
\usepackage{subfigure}
\usepackage{epic}
\usepackage{authblk}
\usepackage{textcomp}
\usepackage[all]{xy}

\usepackage[hypertexnames=false,colorlinks=true,linkcolor=blue,citecolor=blue]{hyperref}
\usepackage[numbers,comma,square,sort&compress]{natbib}
\usepackage[letterpaper,margin=1in,centering]{geometry}

\usepackage[skip=10pt,font=footnotesize]{caption}
\graphicspath{{figures/}}

\newtheorem{Lemma}{Lemma}[section]

\newtheorem{Corollary}[Lemma]{Corollary}

\newtheorem{Definition}[Lemma]{Definition}

\newenvironment{Proof}%
 {\begin{trivlist} \item[]{\bf Proof. }}%
 {\hspace*{\fill}$\rule{.4\baselineskip}{.4\baselineskip}$\end{trivlist}}

 {\begin{trivlist}\item[]\textbf{Acknowledgments.}}{\end{trivlist}}
 
 \newenvironment{Hypothesis}[1]%
  {\begin{trivlist}\item[]{\bf Hypothesis #1 }\em}{\end{trivlist}}

 {\begin{center}\textbf{Abstract}}{\end{center}}

\makeatletter\@addtoreset{figure}{section}\makeatother

\makeatletter \@addtoreset{equation}{section} \makeatother

\newcommand{\R}{\mathbb{R}}

\newcommand{\Z}{\mathbb{Z}}

\newcommand{\cC}{\mathscr{C}}
\newcommand{\F}{\mathcal{F}}

\newcommand{\mc}[1]{\mathcal{#1}}

\newcommand{\ba}{\begin{align}}
\newcommand{\ea}{\end{align}}

\newcommand{\rmd}{\mathrm{d}}

\newcommand{\rme}{\mathrm{e}}
\newcommand{\rmo}{{\scriptstyle\mathcal{O}}}

\newcommand{\eps}{{\varepsilon}}

\newcommand{\uchi}{\underline{\chi}}

\title{Corrigendum to ``Center Manifolds without a Phase Space''}
\author[1]{Gr\'egory Faye\footnote{Corresponding author, gregory.faye@math.univ-toulouse.fr}}
\affil[1]{\small CNRS, UMR 5219, Institut de Math\'ematiques de Toulouse, 31062 Toulouse Cedex, France}
\author[2]{Arnd Scheel\footnote{AS was partially supported by the National Science Foundation through grant  NSF-DMS-1907391.}}
\affil[2]{\small University of Minnesota, School of Mathematics, 206 Church Street S.E., Minneapolis, MN 55455, USA}

\begin{document}

\maketitle

\begin{abstract}
We correct the choice of cut-off function in our construction of center manifolds in  [G. Faye and A. Scheel, \textit{Trans. Amer. Math. Soc.}, 370 (2018), pp. 5843--5885]. The main result there establishes center manifolds for systems of nonlinear functional equations posed on the real line with nonlocal coupling through convolution operators. In the construction, we need to modify the nonlinearity such that it maps spaces of exponentially growing functions into itself and possesses a small Lipschitz constant. The cut-off function presented there is suitable for $L^\infty$-based construction but insufficient for the choice of $H^1$-spaces. We correct this choice and demonstrate the effect on a pointwise superposition nonlinearity. 

\end{abstract}

\section{Construction of a cut-off operator in $H^1_{-\eta}$}

In \cite{FS18}, we aim to carry out a contraction argument to obtain fixed points in spaces of functions $H^1_{-\eta}(\R,\R^n)$, allowing for exponential growth of functions at $\pm\infty$. In general, superposition operators do not map spaces of growing functions into the same space. For example, the map $u(\cdot) \mapsto u(\cdot)^2$, $H^1_{-\eta}(\R,\R)\to H^1_{-\eta}(\R,\R)$ is not defined at $u(x)=\exp(\eta|x|/2)$. One therefore needs to modify the nonlinearity outside of a small ball in the function space where we obtain solutions. Such a cut-off is commonly performed by modifying the nonlinearity outside a small neighborhood of the origin. In \cite{FS18}, we base our cut-off on a \emph{pointwise} modification, changing for instance the quadratic $g:u\mapsto u^2$ as a pointwise evaluation function such that $g(u)=u^2$ for $|u|<\delta$ and $g(u)=0$ for $|u|>2\delta$, for some $\delta$ sufficiently small. This choice however does not guarantee that the Lipschitz constant of $g$, acting as a superposition operator on $H^1$ is small, a property needed in the proof of the main argument there. To see this, choose $u$ a sawtooth compactly supported on  $[0,1]$ with $|u'|=1/\eps$ and $|u|\leq \delta/2$, and $v=u+\delta'$ on the support of $u$. Then, for $\delta'<\delta/2$,  on the support of $u$, 
\[
 \|g(u)-g(v)\|_{H^1}=\|u^2-(u+\delta')^2\|_{H^1}\geq 2\delta' \|u'\|_{L^2}=2\delta'/\eps,
 \]
yielding a lower bound $2/\eps$ on the Lipschitz constant. 

As a consequence, the modification to the nonlinearity in \cite{FS18} does not allow one to carry out the key contraction argument in the proof. The cut-off would be appropriate in an $L^\infty$-based argument, which would in addition require less smoothness in the nonlinearity. On the other hand, the key Fredholm argument for the linear operator from \cite{faye-scheel:13} is only available for $L^2$-based spaces at this point, thus necessitating a different choice of cut-off that we shall introduce next.

We will use a cut-off operator $\chi_\eps$ that cuts off the nonlinearity outside of an $\eps$-ball in $H^1_\mathrm{u}(\R,\R^n)$, the space of functions in $H^1_\mathrm{loc}(\R,\R^n)$ with finite norm 
\[
\|u\|_{H^1_\mathrm{u}}=\sup_y \|u(\cdot -y)\|_{H^1([0,1],\R^n)}.
\]
Note that $H^1_\mathrm{u}(\R,\R^n)\subset H^1_{-\eta}(\R,\R^n)$ for any $\eta>0$.
To construct the cut-off operator, let $\uchi$ be a smooth version of the characteristic  function  of $[-1,1]$, that is 
\[
\uchi(x)=0,\ |x|>2, \qquad \uchi(x)=1,\ |x| <1,\qquad 0\leq |\uchi'(x)|\leq 2.
\]
Next, let $\theta:\R\to\R$ be the smooth generator of a $\Z$-invariant partition of unity on $\R$, that is,
\[
\sum_{j\in\Z} \theta(x-j)=1,\qquad \mathrm{supp}\,\theta\subset \left(-\frac 1 4,\frac 5 4\right
),\qquad \theta(x)\geq 0,\qquad \theta([0,1])\subset \left[\frac 1 2, 1\right].
\]
Note that
\begin{equation}\label{e:cont}
\int_{y\in\R}\theta(x-y)\,\rmd y=\int_{y=0}^1 \sum_{j\in\Z}\theta(x-j-y)\,\rmd y=\int_{y=0}^1 1\,\rmd y=1.
\end{equation}
\begin{Definition}[Cut-off operator]\label{defcutoff}
We define the cut-off operator
\[
\chi:H^1_{-\eta}(\R,\R^n) \to H^1_\mathrm{u}(\R,\R^n)\subset H^1_{-\eta}(\R,\R^n),\qquad  u\mapsto \chi(u),
\]
through 
\[
[\chi(u)](x)=\int_{y\in\R} \uchi\left(\|\theta(\cdot-y)u(\cdot)\|_{H^1(\R,\R^n)}\right)\theta(x-y)u(x)\rmd y.
\]
\end{Definition}
We note that the integral can be understood as a Riemann integral with values in $H^1_\mathrm{u}$ since the integrand is continuous (in fact smooth) in $y$ with values in this space.

We collect some properties of $\chi$. 

\begin{Lemma}[Properties of the cut-off]\label{l:1}
Fix $\eta>0$. There exist    constants $C_0,C_1,C_L>0$ such that
\begin{enumerate}
    \item $\chi$ is well defined as a map from $H^1_{-\eta}(\R,\R^n)$ into $ H^1_{-\eta}(\R,\R^n)$;
    \item $\chi\circ \tau_y=\tau_y\circ \chi$, where $(\tau_yu)(x)=u(x-y)$;
    \item $\chi(u)=u$ if $u\in H^1_\mathrm{u}$ and $\|u\|_{H^1_\mathrm{u}}<C_0$;
    \item $\chi(u)\in H^1_\mathrm{u}$ and $\|\chi(u)\|_{H^1_\mathrm{u}}<C_1$.
    \item $\chi:H^1_{-\eta}(\R,\R^n)\to H^1_{-\eta}(\R,\R^n)$ is globally Lipschitz continuous with Lipschitz constant $C_L$.
\end{enumerate}
\end{Lemma}
\begin{Proof}
To see (i), we estimate 

\begin{equation*}
    \rme^{-\eta |j|}\|\chi(u)\|_{H^1([j,j+1],\R^n)} \leq \int_{y\in[j-2,j+2]}\rme^{-\eta|j|}\|\theta(\cdot-y)u(\cdot)\|_{H^1([j,j+1],\R^n)} 
    \leq 4 \rme^{-\eta|j|} \|u\|_{H^1_\mathrm{u}},
\end{equation*}

using that $|\uchi|\leq 1$, $\sup_\xi|\theta(\xi)|\leq 1$, and that $\mathrm{supp}\,\theta\subset(-\frac 1 2,\frac 3 2)$. 

Next, (ii) is immediate from 
\begin{align*}
[\chi(u)](x-\xi)&=\int_{y\in\R} \uchi\left(\|\theta(\cdot-y)u(\cdot)\|_{H^1(\R,\R^n)}\right)\theta(x-y-\xi)u(x-\xi)\rmd y\\
&=\int_{\tilde{y}\in\R}\uchi\left(\|\theta(\cdot-\tilde{y}+\xi)u(\cdot)\|_{H^1(\R,\R^n)}\right)\theta(x-\tilde{y})u(x-\xi)\rmd \tilde{y}\\
&=\int_{\tilde{y}\in\R}\uchi\left(\|\theta(\cdot-\tilde{y})u(\cdot-\xi)\|_{H^1(\R,\R^n)}\right)\theta(x-\tilde{y})u(x-\xi)\rmd \tilde{y}\\
&=\chi(u(\cdot-\xi))(x),
\end{align*}
where we used substitution $\tilde{y}=y+\xi$ in the second equality and translation invariance of the $H^1$-norm in the third equality. 

Property (iii) follows after exploiting the fact that $
\|\theta(\cdot-y)u(\cdot)\|_{H^1(\R,\R^n)}\|\leq 3 \|u\|_{H^1_\mathrm{u}}<1$ if we choose $C_0\leq 1/3$, hence 
\[
\uchi\left(\|\theta(\cdot-y)u(\cdot)\|_{H^1(\R,\R^n)}\right)=1, 
\]
which results in 
\begin{align*}
\chi(u)&=\int_{y\in\R} \uchi\left(\|\theta(\cdot-y)u(\cdot)\|_{H^1(\R,\R^n)}\right)\theta(x-y)u(x)\rmd y\\
&=\int_{y\in\R} \theta(x-y)u(x)\rmd y=u(x),
\end{align*}
where we used \eqref{e:cont} in the last identity. 

Similarly, 
\begin{align*}
\|\chi(u)\|_{H^1([\xi,\xi+1],\R^n)}&\leq \int_{y\in\R} \uchi\left(\|\theta(\cdot-y)u(\cdot)\|_{H^1(\R,\R^n)}\right)\|\theta(\cdot-y)u(\cdot)\|_{H^1([\xi,\xi+1],\R^n)}\rmd y\\
&\leq \int_{y\in[\xi-2,\xi+2]} \uchi\left(\|\theta(\cdot-y)u(\cdot)\|_{H^1(\R,\R^n)}\right)\|\theta(\cdot-y)u(\cdot)\|_{H^1(\R,\R^n)}\rmd y\\
&\leq 8,
\end{align*}
where the last inequality follows from the fact that $\chi(v)v\leq 2$.

The last estimate (v) is crucial and requires a bit more work. We use that \[
\|\chi(u)\|^2_{H^1_{-\eta}(\R,\R^n)} \lesssim \sum_{j\in\Z} \rme^{-2\eta|j|}\|\chi(u)\|^2_{H^1([j,j+1],\R^n)}\]
and estimate
\begin{align*}
    \rme^{-\eta|j|}&\|\chi(u)-\chi(v)\|_{H^1([j,j+1],\R^n)}\\
    &\leq \int_{y\in[j-2,j+2]} \rme^{-\eta|j|}
    \left\|
    \uchi\left(\|\theta(\cdot-y)u(\cdot)\|_{H^1}\right)\theta(\cdot-y)u(\cdot) -
    \uchi\left(\|\theta(\cdot-y)v(\cdot)\|_{H^1}\right)\theta(\cdot-y)v(\cdot) 
    \right\|_{H^1([j,j+1],\R^n)}\rmd y.
\end{align*}
We estimate the integrand for fixed $y$, uniformly in $y$, which will yield the desired result. 
We may assume that $\|\theta(\cdot-y)v(\cdot)\|_{H^1(\R,\R^n)}<2$ since the integrand vanishes when both arguments of $\uchi$ are larger than $2$, and possibly swapping the role of $u$ and $v$. We further split
\begin{align*}
   \left\|
\uchi\left(\|\theta(\cdot-y)\right.\right.&\left.\left.u(\cdot)\|_{H^1}\right)\theta(\cdot-y)u(\cdot) -    \uchi\left(\|\theta(\cdot-y)v(\cdot)\|_{H^1}\right)\theta(\cdot-y)v(\cdot) 
    \right\|_{H^1([j,j+1],\R^n)}\\
    &\leq 
   \left\|
\uchi\left(\|\theta(\cdot-y)u(\cdot)\|_{H^1}\right)\theta(\cdot-y)u(\cdot) -    \uchi\left(\|\theta(\cdot-y)u(\cdot)\|_{H^1}\right)\theta(\cdot-y)v(\cdot) 
    \right\|_{H^1([j,j+1],\R^n)}\\
    &+\left\|
\uchi\left(\|\theta(\cdot-y)u(\cdot)\|_{H^1}\right)\theta(\cdot-y)v(\cdot) -    \uchi\left(\|\theta(\cdot-y)v(\cdot)\|_{H^1}\right)\theta(\cdot-y)v(\cdot) 
    \right\|_{H^1([j,j+1],\R^n)}\\
    =& \quad \mathrm{I}+\mathrm{II}.
\end{align*}
The first term is easily estimated by
\[
\mathrm{I}\leq \|u-v\|_{H^1([j-2,j+2])}, 
\]
which, together with the exponential weight gives the desired Lipschitz estimate on this part. 

The second term, using a Lipschitz estimate on norms and $\uchi$ yields
\begin{align}
\mathrm{II}&\leq \left|\uchi\left(\|\theta(\cdot-y)u(\cdot)\|_{H^1}\right)-\uchi\left(\|\theta(\cdot-y)v(\cdot)\|_{H^1}\right)\right| \|\theta(\cdot -y)v(\cdot)\|_{H^1([j-2,j+2])}\nonumber\\
&\leq C \left|\|\theta(\cdot-y)u(\cdot)\|_{H^1}-\|\theta(\cdot-y)v(\cdot)\|_{H^1}\right| \|\theta(\cdot -y)v(\cdot)\|_{H^1([j-2,j+2])}\nonumber\\
&\leq C \|u-v\|_{H^1([j-2,j+2])}\|\theta(\cdot -y)v(\cdot)\|_{H^1([j-2,j+2])}\label{e:II}
\end{align}
as $y\in[j-2,j+2]$.

Multiplying by the exponential weight and using the fact that $\|\theta(\cdot-y)v(\cdot)\|_{H^1(\R,\R^n)}<2$, this again gives the desired estimate. 
\end{Proof}

\begin{Corollary}[Localized cutoff]\label{c:1}
The scaled  cutoff 
\[
\chi_\eps(u):=\eps\chi(u/\eps)
\]
satisfies properties. More precisely, fix $\eta>0$. There exist    universal constants $C_0,C_1,C_L>0$ such that
\begin{enumerate}
    \item $\chi_\eps$ is well defined as a map from $H^1_{-\eta}(\R,\R^n)$ into $ H^1_{-\eta}(\R,\R^n)$;
    \item $\chi_\eps\circ \tau_y=\tau_y\circ \chi_\eps$, where $(\tau_yu)(x)=u(x-y)$;
    \item $\chi_\eps(u)=u$ if $u\in H^1_\mathrm{u}$ and $\|u\|_{H^1_\mathrm{u}}<C_0\eps$;
    \item $\chi_\eps(u)\in H^1_\mathrm{u}$ and $\|\chi(u)\|_{H^1_\mathrm{u}}<C_1\eps$.
    \item $\chi_\eps:H^1_{-\eta}(\R,\R^n)\to H^1_{-\eta}(\R,\R^n)$ is globally Lipschitz continuous with Lipschitz constant $C_L$,
\end{enumerate}
\end{Corollary}

\begin{Proof} Most of the properties are immediately clear. The bound on the Lipschitz constant follows immediately from the fact that the Lipschitz constant of compositions is bounded by the product of Lipschitz constants, such that conjugation of a map by scaling does not increase the Lipschitz constant. 
\end{Proof}

\section{Necessary changes to hypotheses in \cite{FS18}}

We may now define the modified nonlinearity from \cite{FS18} as 
\[
\mc F^\eps(u):=\mc F(\chi_\eps(u)),
\]
with $\chi_\eps$ from Corollary~\ref{c:1} in place of (2.4) from  \cite{FS18}.  Hypothesis (H2) in \cite{FS18} is modified by keeping (ii) and (iii) and assuming in addition (3.6)(b).  Note that this latter estimate can usually be verified by exploiting that $\mathcal{F}$ is Lipschitz continuous on $H^1_\mathrm{u}$ in the topology of $H^1_{-\eta}$, that is, 
\[
\|\F(u)-\F(v)\|_{H^1_{-\eta}}\leq L(\|u\|_{H^1_\mathrm{u}},\|v\|_{H^1_\mathrm{u}})\|u-v\|_{H^1_{-\eta}}.
\]
Moreover, $L(a,b)\leq C(a+b)$ if the derivative of $\mc F$ vanishes. Therefore, composing $\F$ with $\chi_\eps$, we find a Lipschitz constant of order $\eps$. 

To be specific, we replace Hypothesis (H2) in \cite{FS18} with the following assumptions. 

\begin{Hypothesis}{(H2).}
We assume that there exists $k\geq 2$ and $\eta_0>0$ such that for all $\epsilon>0$, sufficiently small, the following properties hold.
\begin{enumerate}
\item $\F^\eps$ commutes with translations, $\F^\eps\circ\tau_\xi=\tau_\xi\circ\F^\eps$ for all $\xi\in\R$;
\item $\F^\epsilon:  H^1_{-\zeta}(\R,\R^n) \longrightarrow H^{1}_{-\eta }(\R,\R^n)$ is $\cC^k$ for all nonnegative pairs $(\zeta,\eta)$ such that $0<k\zeta <\eta<\eta_0$,  $D^j_u\F^\epsilon(u): (H^1_{-\zeta}(\R,\R^n))^j \longrightarrow H^{1}_{-\eta }(\R,\R^n)$ is bounded for $0<j\zeta\leq \eta<\eta_0$, $0\leq j\leq k$ and Lipschitz in $u$ for $1\leq j \leq  k-1$;
\item $\F^\epsilon(0)=0$, $D_u\F^\epsilon(0)= 0$ and, as $\epsilon \to 0$, one has the estimate
\begin{equation}
\delta_1(\epsilon):=\textnormal{Lip}_{H^{1}_{-\eta }(\R,\R^n)}(\F^\epsilon)=\mathcal{O}(\epsilon). \label{estLip}
\end{equation}
\end{enumerate}
\end{Hypothesis}
Note that, as a consequence of the Lipschitz bound \eqref{estLip} and the fact that $\F^\epsilon(0)=0$, we readily obtain
\begin{equation}
\delta_0(\epsilon) := \underset{u \in H^1_{-\eta}(\R,\R^n)}{\sup}\| \F^\epsilon(u) \|_{H^{1}_{-\eta }(\R,\R^n)}= \mathcal{O}(\epsilon^2). 
\label{estF}
\end{equation}

Then, directly from our modified Hypothesis (H2), we deduce that the map $\mathcal{S}^\epsilon(u,u_0)$ from (3.5) in \cite{FS18} satisfies the estimates
\begin{align*}
\|\mathcal{S}^\epsilon(u,u_0)\|_{H^1_{-\eta}(\R,\R^n)} &\leq C(\eta) \left( \delta_0(\epsilon)+\|u_0\|_{H^1_{-\eta}(\R,\R^n)}\right),\\
\|\mathcal{S}^\epsilon(u,u_0)-\mathcal{S}^\epsilon(v,u_0)\|_{H^1_{-\eta}(\R,\R^n)} &\leq C(\eta) \delta_1(\epsilon)\|u-v\|_{H^1_{-\eta}(\R,\R^n)},
\end{align*}
for all $u,v\in H^1_{-\eta}(\R,\R^n)$ and $u_0\in\mathcal{E}_0$.

Hypothesis (H2$\mu$) also needs to be adapted according to the previous modifications, defining the  parameter dependent nonlinearity through
\[
\mc F^\eps(u,\mu):=\mc F(\chi_\eps(u),\mu).
\]

\begin{Hypothesis}{(H2$\mu$).}
We assume that there exists $k\geq 2$ and $\eta_0>0$ such that for all $\epsilon>0$, sufficiently small, the following properties hold.
\begin{enumerate}
\item $\F^\eps(\cdot,\mu)$ commutes with translations for all $\mu$, $\F^\eps\circ\tau_\xi=\tau_\xi\circ\F^\eps$ for all $\xi\in\R$;
\item $\F^\epsilon:  H^1_{-\zeta}(\R,\R^n)\times \mathcal{V}_\mu \longrightarrow H^{1}_{-\eta }(\R,\R^n)$ is $\cC^k$ for some $0\in \mathcal{V}_\mu\subset \R^p$ and for all nonnegative pairs $(\zeta,\eta)$ such that $0<k\zeta <\eta<\eta_0$,  $D^j_u\F^\epsilon(u,\mu): (H^1_{-\zeta}(\R,\R^n))^j \longrightarrow H^{1}_{-\eta }(\R,\R^n)$ is bounded for $0<j\zeta\leq \eta<\eta_0$, $0\leq j\leq k$ and Lipschitz in $u$ for $1\leq j \leq  k-1$ uniformly in $\mu\in\mathcal{V}_\mu$;
\item $\F^\epsilon(0,0)=0$, $D_u\F^\epsilon(0,0)= 0$ and, as $\epsilon \to 0$ one has the estimates
\begin{equation}
\delta_1(\epsilon):=\textnormal{Lip}_{H^{1}_{-\eta }(\R,\R^n)\times  \mathcal{V}_\mu}(\F^\epsilon)=\mathcal{O}(\epsilon +\|\mu\|). \label{estLipmu}
\end{equation}
\end{enumerate}
\end{Hypothesis}

\section{Verifying Hypothesis (H2) for the superposition operator $u\mapsto u^2$}

We demonstrate that the new hypotheses can be verified for the simplest case of the quadratic  superposition operator $\F:u\mapsto u^2$ which leads to $\F^\eps(u):=\chi_\eps(u)^2$. Adaptations for more general pointwise superposition operators are straightforward, as are compositions with convolution operators.

Property  (i) follows from the translation invariance property of the cutoff. Regarding (iii), clearly $\F^\eps(0)=0$. We will discuss the vanishing derivative after establishing (ii). The estimates in (iii) can readily be verified as follows.  We have that
\[ \| \chi_\eps(u)^2 \|_{H^1_{-\eta}}\lesssim \| \chi_\eps(u)^2 \|_{H^1_\mathrm{u}}\lesssim \| \chi_\eps(u) \|_{H^1_u}^2 \lesssim \epsilon^2,
\]
using that $\chi_\eps$ maps into an $\eps$-ball in $H^1_\mathrm{u}$, which is an algebra. 
For the Lipschitz constant, we use that 
 $|uv|_{H^1_{-\eta}}\leq |u|_{H^1_\mathrm{u}}|v|_{H^1_{-\eta}}$:
\[
\|\chi_\eps(u)^2-\chi_\eps(v)^2 \|_{H^1_{-\eta}} \leq \|\chi_\eps(u)+\chi_\eps(v) \|_{H^1_\mathrm{u}} \|\chi_\eps(u)-\chi_\eps(v) \|_{H^1_{-\eta}} \lesssim \epsilon \|\chi_\eps(u)-\chi_\eps(v) \|_{H^1_{-\eta}}.
\]
We now establish differentiability as stated in (ii), which roughly follows the usual arguments but is somewhat lengthy because of the cumbersome expression for the cutoff operator. Let $0<\zeta<\eta$ and prove that $\F^\epsilon:  H^1_{-\zeta}(\R,\R^n) \longrightarrow H^{1}_{-\eta }(\R,\R^n)$ is $\cC^1$. Higher differentiability can then be established in an analogous fashion. It is sufficient to prove differentiability for $\eps=1$, given that the scaling operators are linear and differentiable on any of the exponentially weighted spaces. We define the candidate for the derivative of 
\[
L(u)\cdot v = 2\chi(u) \left(\chi^1(u)\cdot v\right),
\]
where
\begin{align*}
\left[\chi^1(u)\cdot v\right](x)&:=\left(\int_{y\in\R} \uchi'\left(\rho_y(u)\right)D_u\rho_y(u)v \,\theta(x-y)u(x)\rmd y\right)\\
&+\int_{y\in\R} \uchi\left(\rho_y(u)\right)\theta(x-y)v(x)\rmd y,
\end{align*}
and we have set 
$$\rho_y(u):=\|\theta(\cdot-y)u(\cdot)\|_{H^1(\R,\R^n)},$$ and 
\[
D_u\rho_y(u)v:=\left\langle \frac{\theta(\cdot-y)u(\cdot)}{\rho_y(u)},\theta(\cdot-y)v(\cdot) \right\rangle.
\] 
We claim that
\begin{itemize}
    \item[(A)] $L:H^1_{-\zeta}(\R,\R^n) \longrightarrow \mathcal{B}\left(H^1_{-\zeta}(\R,\R^n),H^{1}_{-\eta }(\R,\R^n)\right)$ is continuous;
    \item[(B)] $L(u)v$ is the Gateaux derivative of $\F: H^1_{-\zeta}(\R,\R^n) \longrightarrow H^{1}_{-\eta }(\R,\R^n)$ in the direction of $v$.
\end{itemize}
Both statements together imply Fr\'echet differentiability of $\F$.

To see (A), a direct calculation exploiting that $\uchi'=0$ unless $\rho_y\in (1,2)$ shows that 
\[L(u)\in  \mathcal{B}\left(H^1_{-\zeta}(\R,\R^n),H^{1}_{-\zeta }(\R,\R^n)\right)\subset \mathcal{B}\left(H^1_{-\zeta}(\R,\R^n),H^{1}_{-\eta }(\R,\R^n)\right).\]

To prove continuity of $L$, we split
\[
\|L(u_1)v-L(u_2)v\|_{H^1_{-\eta}}\lesssim \|(\chi(u_1)-\chi(u_2))\chi^1(u_1)v)\|_{H^1_{-\eta}}+ \|\chi(u_2)(\chi^1(u_1)-\chi^1(u_2))v\|_{H^1_{-\eta}}=:\mathrm{I}+\mathrm{II}. 
\]
The first term can be estimated as follows. For $\delta=\eta-\zeta>0$ small, we note that
\[
\|(\chi(u_1)-\chi(u_2))\chi^1(u_1)v)\|_{H^1_{-\eta}}\leq \|(\chi(u_1)-\chi(u_2))\|_{H^1_{-\delta}} \|\chi^1(u_1)) v\|_{H^1_{-\zeta}}.
\]
 Let $\beta\in(0,1)$ such that $1-\frac{\beta\delta}{\zeta}>0$. Using the shorthand notation $H^1_j=H^1([j-2,j+2],\R^n)$, we can estimate
\begin{align*}
    \|\chi(u_1)-\chi(u_2)\|^2_{H^1_{-\delta}} &\lesssim  \sum_j \left(\rme^{-\delta|j|} \|\chi(u_1)-\chi(u_2)\|_{H^1_j}\right)^2\\
      = & \sum_j \left(\rme^{-2\zeta|j|}  \|\chi(u_1)-\chi(u_2)\|_{H^1_j}^2\right)^{\frac{\beta\delta}{\zeta}}\rme^{-2(1-\beta)\delta|j|}  \left(\|\chi(u_1)-\chi(u_2)\|_{H^1_j}^2\right)^{1-\frac{\beta\delta}{\zeta}}.
\end{align*}
Next, using the boundedness of $\chi$, we readily see that
\[
\left(\|\chi(u_1)-\chi(u_2)\|_{H^1_j}^2\right)^{1-\frac{\beta\delta}{\zeta}} \lesssim \left(\|\chi(u_1)-\chi(u_2)\|_{H^1_{ul}}^2\right)^{1-\frac{\beta\delta}{\zeta}} \leq 2^{2\left(1-\frac{\beta\delta}{\zeta} \right)}.
\]
As a consequence, we obtain 
\begin{align*}
    \|\chi(u_1)-\chi(u_2)\|^2_{H^1_{-\delta}} &\lesssim \sum_j \left(\rme^{-2\zeta|j|}  \|\chi(u_1)-\chi(u_2)\|_{H^1_j}^2\right)^{\frac{\beta\delta}{\zeta}}\rme^{-2(1-\beta)\delta|j|}  \\
   & \lesssim \left( \sum_j \left(\rme^{-\zeta|j|}  \|\chi(u_1)-\chi(u_2)\|_{H^1_j}\right)^2 \right)^{\frac{\beta\delta}{\zeta}} \left( \sum_j \rme^{-2 \frac{(1-\beta)\delta \zeta}{\zeta-\beta\delta}|j|}\right)^{1-\frac{\beta\delta}{\zeta}}\\
   & \lesssim  \|\chi(u_1)-\chi(u_2)\|_{H^1_{-\zeta}}^{2\frac{\beta\delta}{\zeta}} \\&\lesssim  \|u_1-u_2\|_{H^1_{-\zeta}}^{2\frac{\beta\delta}{\zeta}},
\end{align*}
where we have used  H\"older's inequality and  Lipschitz continuity of $\chi$ in $H^1_{-\zeta}$. Finally, we notice that $\chi^1(u_1)$ is a bounded linear operator 
\[
\|\chi^1(u_1) v\|_{H^1_{-\zeta}} \leq \|\chi^1(u_1)\|_{\mathcal{B}(H^1_{-\zeta},H^1_{-\zeta})} \|v\|_{H^1_{-\zeta}}
\]
for all $u_1\in H^1_{-\zeta}$ as from its definition  we have
\begin{align*}
\|\chi^1(u_1)) v\|_{H^1_{-\zeta}}^2&\lesssim \sum_j\rme^{-2\delta|j|} \|\chi^1(u_1)) v \|_{H^1([j,j+1],\R^n)}^2\\
&\lesssim \sum_j\rme^{-2\delta|j|} \|v \|_{H^1([j-3,j+3],\R^n)}^2.
\end{align*}
 This further implies that $\|\chi^1(u_1)\|_{\mathcal{B}(H^1_{-\zeta},H^1_{-\zeta})}$ is independent of $u_1\in H^1_{-\zeta}$. We conclude that $(\mathrm{I})$ gives H\"older continuity.
 

The second term (II) can be dealt with in a similar fashion. Inspecting the definition of $\chi^1$, we see that the second summand $\int_{y\in\R} \uchi\left(\rho_y(u)\right)\theta(x-y)v(x)\rmd y$ can be treated in an identical fashion and it will be sufficient to establish continuity of $\uchi'\left(\rho_y(u)\right)D_u\rho_y(u)v$ in $H^1_{-\eta}$, uniformly in  $v\in H^1_{-\zeta}$, bounded. This can be seen again by splitting up the product, 
\begin{align*}
    \|\uchi'\left(\rho_y(u_1)\right)D_u\rho_y(u_1)v-&\uchi'\left(\rho_y(u_2)\right)D_u\rho_y(u_2)v\|_{H^1_{-\eta}} \\
    \leq & \|\uchi'\left(\rho_y(u_1)\right)D_u\rho_y(u_1)v-\uchi'\left(\rho_y(u_2)\right)D_u\rho_y(u_1)v\|_{H^1_{-\eta}}\\
    &+ \|\uchi'\left(\rho_y(u_2)\right)D_u\rho_y(u_1)v-\uchi'\left(\rho_y(u_2)\right)D_u\rho_y(u_2)v\|_{H^1_{-\eta}}\\
    =:&\ \mathrm{(III)}+\mathrm{(IV)}.
\end{align*}

Since $\uchi'$ is supported in $[1,2]$, it is sufficient to verify this when $\rho_y(u)\in[1,2]$, which implies  that $D_u\rho_y(u)v$ is Lipschitz and bounded in $u$ and H\"older continuity of (IV). On the other hand, (III) can be estimated using \eqref{e:II} and the subsequent arguments for (I). 

To see (B), we need to show that  
\[
\|\F(u+\tau v)-\F(u)-\tau D_u\F(u)v\|_{H^1_{-\eta}}=\rmo(\tau) \text{ for } v\in{H^1_{-\zeta}}.
\]

This in turn follows immediately from differentiability of $\uchi$ and the estimates for continuity of $L(u)$ using the fundamental theorem of calculus.  

Lastly, the fact that $D\F^\eps(0)=0$ follows immediately from $L(0)=0$. 

\section*{Acknowledgments} We are grateful to Tien Truong, Erik Wahl\'en and Miles Wheeler who kindly pointed out the inconsistencies with the cut-off function in our original work and provided valuable feedback in the process of writing this corrigendum.

\bibliography{plain}

\end{document}